\documentclass{elsart}
\journal{ \ \ }
\renewcommand{\theequation}{\thesection.\arabic{equation}}
\usepackage[centertags]{amsmath}
\usepackage{newlfont}
\usepackage{color}
\usepackage{graphicx}
\usepackage{psfrag}
\usepackage{float}
\usepackage{caption}
\usepackage[english]{babel}
\usepackage{multicol,caption}
\usepackage{caption}
\usepackage{subcaption}
\usepackage[centertags]{amsmath}
\usepackage{psfrag}
\usepackage{listings}
\usepackage{booktabs}
\usepackage[T1]{fontenc}
\usepackage[english]{babel}
 \usepackage{hyperref}
\usepackage{amsfonts}
\usepackage{xcolor}
\usepackage{hyperref}
\usepackage[utf8]{inputenc}
\usepackage[english]{babel}
 
\newtheorem{theorem}{Theorem}[section]

\floatstyle{plain}
\restylefloat{figure}

\begin{document}
\begin{frontmatter}
\title{Nested Performance Profiles for Benchmarking  Software}

\author{ Rasoul Hekmati$^1$, Hanieh Mirhajianmoghadam $^2$}
\address{1. Department of Mathematics, University of Houston,TX, USA  rhekmati@math.uh.edu  \\ 2. Collage of Optometry, University of Houston}

\begin{abstract}
In order to compare and benchmark the mathematical software, the performance profiles have been introduced \cite{more}. However, it has been proved that the algorithm is not flawless. The main issue with the performance profile is that it may rank the solvers with respect to the best solver, by excluding the best one and running the algorithm on the remaining set of the solvers, the method may rank the solvers in a different way. We characterize such systems of problems-solvers and propose an efficient and reliable algorithm to overcome this negative side effect. The proposed method is unbiased in comparing the solvers and is successful in detecting the top ones. 
\end{abstract}
\end{frontmatter}
%\begin{multicols}{2}
\section{Introduction}

\quad For a set of mathematical software (such as optimization packages), there are several available solvers. Each solver shows a  superior performance on some of the problems and inferior performance on other problems. This makes it difficult to determine which one is better. The interpretation and analysis of the data generated by the benchmarking process have been discussed by Dolan and More \cite{more}. Many benchmarking efforts involve tables displaying the performance of each solver on each problem for a set of metrics such as CPU time, the number of function evaluations and iteration counts for algorithms. The solver's average or cumulative total for each performance metric over all the problems is sometimes used to evaluate performance \cite{benson,bongartz,conn}. As a result, a small number of  difficult problems can influence the overall performance. 
In the 1990s, some researchers ranked the solvers \cite{bongartz,conn,nash}. They counted the number of times that a solver comes
in $k^{th}$ place,  for $k = 1, 2, 3$. Ranking the solvers' performance for each problem helped prevent a minority of  problems from influencing the results. Information on the size of the improvement, however, was lost \cite{more}. Comparing the medians and quartiles of some performance metric has its own disadvantages.
Comparing solvers by the ratio of one solver's performance to the best performance \cite{billups} was also not a flawless approach  (see \cite{more} for detail).
Dolan and More introduced performance profiles (cumulative distribution function of a performance metric) as a tool for evaluating and comparing the performance of mathematical  software. They used the ratio of the computing time of the solver versus the best time of all of the solvers as the performance metric. They showed that performance profiles eliminate the influence of a few of problems on the benchmarking process and the sensitivity of results associated with the ranking of solvers \cite{more}. \\
\indent A negative aspect of performance profiles has been discussed in \cite{gloud}. Gould and Scott showed that if performance profiles are used to compare more than two solvers,  we can determine which solver has the highest probability $\rho_ i ( \tau )$ of
being within a factor $\tau$ of the best solver, but we cannot necessarily determine the performance of one solver relative to another that is not the best. So if we eliminate the best solver and re-do the calculations we may end up with a different answer on the new set of solvers. In some cases, we need to partially rank the solvers. For example, a user might not have access to the best solver and so might want to know which one is the second best solver. Another example is that they might want to identify several top solvers. To overcome this problem, Gould and Scott suggested one option would be to produce a series of performance profiles, excluding the best solver over the range from successive profiles and repeat this procedure until only two remains. This method works if we have a small number of solvers whereas for a large number of solvers this method is not practical (e.g. the problem of selecting the optimal parameter or doing a fine tuning, here each parameter defines a new solver). In this paper, we introduce Nested Performance Profile that combines all the relative features of solvers and gives a single graph to rank the solvers. It uses consecutive performance profiles achieved by eliminating the best solver, the elimination that defines a new reduced system of solvers-problems (which naturally generates nested systems). Nested performance profile is the mean performance over all the reduced systems.
\section{The Method}%%%%%%%%%%%%%%%%%%%%%%%%%%%%
 \quad Consider a set of solvers $S$ on a test set $P$. Let $n_s$ be the number of solvers and $n_p$ the number of problems. We use computing time  as  a performance measure. For each problem $p$ and solver $s$,
$$t_{p,s}:= computing \ time \ required \ to \ solve \ problem \ p \ by \ solver \ s$$
Dolan and More compared the performance on problem $p$ by solver $s$ with the best performance by any solver on this problem, i.e. they used the performance ratio:
\begin{equation}
r_{p,s}=\frac{t_{p,s}}{min\{t_{p,s}:s \in S\}}.
\end{equation}
Assume that a parameter $r_M \geq r _{p,s} $ for all $p, s$ is chosen, and $r_{ p,s} = r_ M$ if and only if solver $s$ does not solve problem $p$. It has been shown that the choice of $r _M$ does not affect the performance evaluation \cite{more}.\\
\indent The probability for solver $s \in S$ that a performance ratio $r_{ p,s}$ is within a factor $\tau \in R$ of the best possible ratio is defined as
\begin{equation}
\rho_s(\tau)=\frac{1}{n_p} | \{p \in P: r_{p,s} \leq \tau \}|.
\end{equation}
The function $\rho_s$ is the cumulative distribution function for the performance ratio. A plot of the performance profile shows the major
performance characteristics. We prefer the solvers with large probability $\rho_s(\tau)$. If we are interested
only in the number of wins, we need only to compare the values of $\rho _s (1)$ for all of the solvers ( $\rho_i (0)$ when a log scaled performance profile has been used ). While $\rho _i ^* = \lim_{\tau \rightarrow \infty } \rho_i ( \tau )$ gives the fraction  for which solver $i$ is successful without considering the speed of convergence. However, if we are interested in solvers with a high probability of success, we should choose those for which $\rho_i^*$ is largest.
Performance profiles are insensitive to the results on a small number of problems, they are also largely unaffected by small changes in results over many problems \cite{more}. \\
\indent  In Nested Performance Profiles in order to rank the top $k$ solvers we have $k$ waves of the performance profiles where  $k<n_s$. In each wave, the best solver is detected and the corresponding performance ratios are saved. Then after eliminating the best solver,  the next wave on the reduced set of solvers is started and  the performance ratios for the eliminated solver(s) are repeated. This defines $k$ performance profiles on $k$ nested sets of solvers. The final performance profile is the mean of the nested profiles. This naturally mitigates the negative side effect of the regular performance profiles, and the achieved graph benchmarks the top $k$ solvers. In this paper, we use the upper index to specify the wave number ,e.g.,  the performance ratios for the $k^{th}$ wave is $r_{p,s}^k$.\\
\indent At the first wave, for simplicity lets set $r_{p,s}=r_{p,s}^1$ for $p\in P, s \in S$ and  $\rho_s=\rho_s^1$ for  $s\in S$. The set of current solvers is $S'=S$ and the set of eliminated solvers is $S^*=\emptyset$. The best solver $s^*$ is identified  by 
 $$s^*=\max_s \{|p\in P : r_{p,s}=1|\},$$ 
or by
$$s^*=\min_s\{\sum_{p\in P} r_{p,s} \}.$$
 The first choice gives the solver with most wins ,i.e., the number of problems for which the solver works best, and the second one gives the solver with best overall mean performance i.e. the mean of ratios. If there is more than one solver with this property then we pick one at random.\\
\indent Now in order to start the second wave, we have to exclude $s^*$ from the set of solvers: $$S'=S' \setminus \{s ^*\}.$$  The updated set of best solvers is  $S^*=S^* \cup s^*$. This naturally defines a new system of solvers and problems (on $S'$), the second wave of performance ratios on the new system launches similar to the first wave:
$$r_{p,s}^2=\frac{t_{p,s}}{min\{t_{p,s}:s \in S\}}  \ \ and\ \  \rho_s^2(\tau)=\frac{1}{p_n} size \{p \in P: r_{p,s}^2 \leq \tau \}. $$
Now  for the eliminated  solver $s^*$ and for a specific problem$p$, if it still shows the best performance (if $r_{p,s^*}^{1}=1$)  we repeat the previous performance ratios:
$$r_{p,s^*}^2=r_{p,s^*}^1,  \ \ \ \ \ s^*\in S^*.$$
\indent If the eliminated solver is not the best solver for a specific problem (if $r_{p,s}^{1}\neq1 $) then the algorithm deals with it like a non-eliminated  solver.
We have to repeat this procedure $k$ times if the top $k$ solvers are  what we are going to specify. Clearly, $k=size(S)-1$ is the value for which the algorithm compares all the solvers. The overall performance profile is the mean of nested performance profiles:
\begin{equation}
\rho_s^{Overall} =\frac{ \sum_{i=1}^k \rho_s^i(\tau)}{k}  , i=1 ... k
\end{equation}
In this way, the comparison is not based on the best solver but it is based on the top $k$ solvers. \\\\
\textbf{Nested Performance Profile Algorithm:}\\
Step 0. Set: the set of best solvers $S^*= \emptyset $ , the set of remaining solvers $S'=S$\\
Step 1. Calculate $r_{p,s}^1$ using (1)\\
Step 2. Calculate $\rho_s^1$ using (2)\\
Step  3. For $i=2,...,k$\\
 \hspace*{1.5ex}3.1 Find the best solver $s^*$ and update  $S^*=S^* \cup s^*$ and  $S'=S'\setminus \{ s ^*\}$ \\
 \hspace*{1ex} 3.2 For $p \in P$ and $s\in S$, repeat\\
\hspace*{4.5ex} \{\\
\hspace*{4.5ex} if $s\in S'$, calculate: \\
\hspace*{7.5ex}$r_{p,s}^i=\frac{t_{p,s}}{min\{t_{p,s}:s \in S\}}  \ \ and\ \  \rho_s^i(\tau)=\frac{1}{p_n} size \{p \in P: r_{p,s}^i \leq \tau \} $\\
\hspace*{4.5ex} else ( if $s\in S^*$):\\
\hspace*{7.5ex} if $r_{p,s}^{i-1}\neq1$:\\
\hspace*{10.5ex}$r_{p,s}^i=\frac{t_{p,s}}{min\{t_{p,s}:s \in S\}}  \ \ and\ \  \rho_s^i(\tau)=\frac{1}{p_n} size \{p\in P: r_{p,s}^i \leq \tau \} $
\hspace*{8ex}else ( if $r_{p,s}^{i-1}=1 $):\\ 
\hspace*{10.5ex}set $r_{p,s}^i=1  $\\
\hspace*{4.5ex} \}\\
Step 4. Calculate the overall performance profile using (3).

\begin{theorem}
The nested performance profiles are insensitive to the results on a small number of problems i.e. if $n_p$ is reasonably large,  then the result on a particular problem $q$ does not greatly affect the nested performance profiles.
\end{theorem} 

{\it Proof.}   As in \cite{more}, if the observed time sets are $t_{p,s}$ and $ \hat{t}_{p,s}$, where $$\hat{t}_{p,s}=t_{p,s}, \ \ \ \ \ \ p \in P \setminus \{q\},$$
for some problem $q \in P$, then $\hat{r}_{p,s}=r_{p,s}$ for $p \in P \setminus  \{q\}$ and for $s \in S$ we have:
$$|\rho_s^i(\tau)-\hat{\rho_s}^i(\tau)| \leq \frac{1}{n_p}    \ \ \ \ \ \ \ \ \ \  i=1, . . . , k , \tau \in \mathbb{R}.$$
 For the overall performance profile:
$$|\rho_s^{Overall} (\tau)-\hat{\rho_s}^{Overall} (\tau)|=\frac{1}{k}|\sum_{i=1}^k \rho_s^i(\tau)- \sum_{i=1}^k \hat{\rho_s}^i(\tau)|\leq\frac{k}{kn_p}=\frac{1}{n_p},$$
moreover $\hat{\rho_s}^{Overall}(\tau) = \rho_s^{Overall}(\tau)$ for $\tau<min\{r_{q,s}^i, \hat{r}_{q,s}^i\}$ or $\tau>max\{r_{q,s}^i, \hat{r}_{q,s}^i\}.$
Thus, if $ n_ p$ is large enough, then the result on a particular problem $q$ does not affect
the nested performance profiles.
\qed
Lets define $Rank_S$ to be the sequence showing the index of the ranked solvers in $S$. Here, $s^i$ is the solver number $i$ and $$t_{P,S^{i,j}} := t_{p,s} \ \ for  \   p \in \ P \ ,  \ s\in \{s^i,s^j\}.$$ 
\begin{theorem}
The performance profiles are sensitive to the elimination of the best solver i.e. if $s^* \in S$ is the best solver then   $Rank_{S \setminus \{s^*\}}$ is not necessarily equal to $ s^* \cup Rank_S \setminus \{s^*\}$. The nested performance profils are not.
\end{theorem} 

{\it Proof.}  Assume that $n_s=3$ and $P_1$, $P_2$ and $P_3$ are  partitions of P, so that $P_1 \cup P_2 \cup P_3=P$, $|P_1|>n_p/2 $, $|P_2|>n/4$   such that $ t_{p,s^1} < t_{p_1,s^{2,3}} $,  $ t_{p_3,s^3} < t_{p_3,s^{1,2}} $ and  $ t_{p_2,s^2} < t_{p_2,s^{1,3}} $ with the extra condition that  $ t_{p_1,s^3} < t_{p_1,s^2} $. Clearly  $r_{P_i,s^i}=1$ for $ i=1,2,3$ considering the size of each partition $\rho _{p,s^i}(1)=|P_i|, i=1,2,3$ i.e. $s^1$ is better than $s^2$ and $s^2$  better than $s^3$  (or $Rank_S=[1,2,3]$). After eliminating the best solver  $s^1$, we have the new system with 2 solvers. We have $\rho_{p,s^3}=|P_1|+|P_3| $ while $\rho_{p,s^2}=|P_2|$ thus $s^3$ is better than $s^2$ i.e. $Rank_{S\setminus \{s^*\}}=[3,2]$. The defined system of solvers and problems proves the theorem. 
For nested performance profile we have: $\rho_{s^i}^1(1) =|P_i|$ and $\rho_{s^2}^2=|p_2|$, $\rho_{s^1}^2=|p_1|$, $\rho_{s^3}^2=|p_1|+|p_3|$ so  $\rho_{s^1}^{Overall}(1)=|P_1|$, $\rho_{s^2}^{Overall}(1)=|p_2|$ and $\rho_{s^3}^{Overall}(1) = |p_3|+\frac{|p_1|}{2}$. Since $|p_1|>=2|p_2|$ so $s^1$ is better than $s^3$ and $s^3$  better than $s^2$.

\indent For the general case with n partitions and n solvers let $P= \bigcup _{i=1}^n P_i$ and $|P_i|>n_p/(2^n)$ and $t_{p_i,s^i}<t_{p_i,S^{1,...,n \setminus \{i\}}}$ with the extra condition that $t_{p_i,s^{i+2}}<t_{p_i,s^{i+1}}, i=1,..., n-2$   by a similar discussion  we can build a system of problems such that the performance profile may rank the solvers in a wrong way.
So the performance profiles are insensitive to changes in results on a small number of problems and sensitive to changes in the set of solvers. They are also largely unaffected by small changes in results over many problems.
\qed

\begin{theorem}
Let $r_i$ and $\hat{r}_i$ for $1 \leq i \leq n_p$ be performance ratios for some solver. Let $\rho$ and $\hat{\rho}$ be, respectively, the nested performance profiles defined by these ratios. If $|r_i - \hat{r}_i | \leq \epsilon$  for some $\epsilon > 0$, then $$\int_1^{\infty}|\rho^{Overall}(t)-\hat{\rho}^{Overall}(t)| dt \leq \epsilon.$$
\end{theorem} 

{\it Proof.}  More and Dolan in \cite{more} proved that the theorem holds for performance profiles, i.e. for each single performance profile: $$\int_1^{\infty}|\rho^i(t)-\hat{\rho}^i(t)| dt \leq \epsilon  , \ \ \ \ i=1,...,k.$$
So,
 $$\int_1^{\infty}|\rho^{Overall}(t)-\hat{\rho}^{Overall}(t)| dt =\frac{1}{k}\int_1^{\infty}|\sum_{i=1}^k(\rho^i(t)-\hat{\rho}^i(t))| dt \leq \frac{1}{k}\int_1^{\infty}\sum_{i=1}^k|\rho^i(t)-\hat{\rho}^i(t)| dt $$
$$=\frac{1}{k}\sum_{i=1}^k \int_1^{\infty}|\rho^i(t)-\hat{\rho}^i(t)| dt \leq \frac{1}{k} k \epsilon = \epsilon.$$
\qed

\section{Numerical Experiments}

\quad The artificial sample data proposed in \cite{gloud} is given in Table 1. Using this data for five test problems and three solvers  and the corresponding logarithmic scaled performance profiles given in Figure 1, we can see the weakness of this method. This system of problems and solvers is what we characterized in the proof of theorem 2.
\begin {table}[h!]
\centering
\begin{tabular}{ |c|c|c|c| } 
 \hline
Problem & Solver A &Solver B  &Solver C    \\ 
\hline
 1&  2&  1.5 &   1 \\ 

 2&1  &  1.2& 2 \\ 

 3& 1 &4   & 2  \\ 

 4& 1 & 5  & 20  \\ 

 5&  2& 5  & 20  \\ 
\hline
\end{tabular}
\caption{The Artificial test set, the smaller the statistics, the better the solver performance}
\label{Table:1}
\end {table}
With $S_1$ = \{Solver $A$, Solver $B$, Solver $C$\}, Solver $A$ is the best on 80\% of the problems, Solver $B$ is not the winner in $\tau \in [0,2]$,
If we are interested in having a solver that can solve at least 60\% of the test
problems with the greatest efficiency, then we should choose solver $A$ or $C$. However,
if $S_2$ = \{Solver $B$, Solver $C$\} (i.e., Solver $A$ is removed), Solver $B$, which was the second
best solver in $S_1$ on 60\% of the test set, is the best solver in $S_2$ \cite{gloud}.
\begin{figure}[h]
\caption{ Performance Profile for the artificial test set}
\centering
\includegraphics[width=3in]{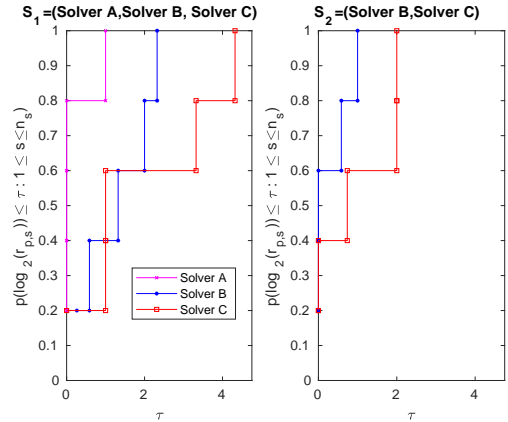}
\end{figure}
The matlab solver (\cite{hek2}) which we used is a modification of the regular performance profile solver written by Dolan and More \cite{cops}. After running the nested performance profile on the sample data (Fig 2), we  successfully ranked the solvers. Solver $A$ is the best solver, solver $B$ is the second one and finally solver $C$ is the last choice. Clearly solver $B$ is superior to Solver $C$ and we don't need to eliminate Solver $A$ to investigate this issue and the nested performance profile could eliminate the relative comparing effect of the regular performance profile on the artificial data. 
\begin{figure}[h]
\caption{ Regular and Nested Performance Profile }
\centering
\includegraphics[width=3in]{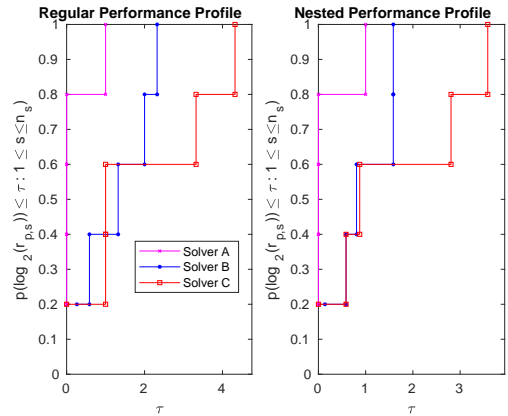}
\end{figure}
As a real example, we used table V in \cite{scott} which demonstrates the Total Time Required for Subset CUTEst problems by each method. The winner is MA87 and the probability that MA87  is the winner on a given
problem is about 61\%. It is noteworthy that we scaled the x-axis and took $0<\tau<0.6* ( max \ Ratio)$ as the nested performance profile has a bigger maximum ratio comparing to the regular performance profile. If we choose to be within a factor of 4 of the best solver, then MA87 is still the best choice; but the performance
profile shows that the probability that this solver can solve the problems within a factor
>5 of the best solver is only about 80\%. Solver "diagonal" has a lower number of wins than "MA87",
 but its performance becomes much more competitive if we extend the $\tau$ of interest to more than 5. "MI35" is the next appropriate solver. The question is: what will happen if we discard the best solver "MA87"? is the solver "diagonal" better than the solver "MI35"? By looking at the nested performance profile (Fig 3. Right) clearly "diagonal" is better than "MI35", something that can not be declared directly by looking at the performance profiles especially for $2<\tau<6 $. \\
\indent Also by looking at the regular performance profile, we may wrongly conclude that "MI35" is better than "None" in solving the problems in a high $\tau$  or we may wrongly conclude that MA87 has the same performance as MIQR, while the regular performance profile ONLY says something about the best solver and no conclusion can be made on the next best solvers.
Clearly, we can rank the solvers by using the nested performance profile. For $\tau>8$, we can observe that the results in nested performance profile are reliable and we don't need to eliminate the best solver and run a sequence of performance profiles to figure out the top solvers.
\begin{figure}[h]
\caption{ Time for Subset CUTEst Problems }
\centering
\includegraphics[width=4in]{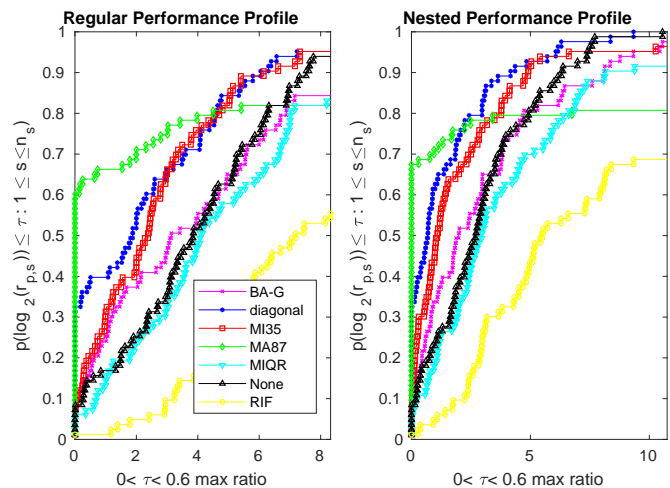}
\end{figure}
\section{Conclusion}
 \quad Performance profile provides a measure to compare multiple solvers. For a binary comparison it is a strong tool  for a selected range of $\tau$. However, if performance profile is used to compare multiple solvers we can determine which one has a higher probability of being within a factor $\tau$ of the best solver but we can not evaluate the performance of one solver relative to another one that is not the best. To address this problem we introduced the nested performance profile that uses consecutive performance profiles achieved by eliminating the best solver and calculates the mean performance over all of the runs. This algorithm combines the relative features of the solvers and gives a reliable criteria to compare all the solvers together. This can be useful if we don't have access to the first solver or if we are interested in determining second or third best solvers out of a large set of solvers. The proposed method is a practical approach to deal with the fine tuning problem which can be seen as a benchmarking problem.
\section{Acknowledgement}
We thank Mozahid Haque for comments that improved the manuscript.

\end{document}